\documentclass[a4paper,12pt]{article}
\usepackage{amsmath}
\usepackage{amsfonts}
\usepackage{amssymb}
\usepackage{lscape}
\usepackage{geometry}
\usepackage[doublespacing]{setspace}
\usepackage{natbib}
\usepackage[hiresbb]{graphicx}
\usepackage{multicol}
\usepackage{float}
\usepackage{quoting}

\begin{document}

\title{New Proofs of the Basel Problem \\
	 using Stochastic Processes}
\author{Uwe Hassler\thanks{\textbf{Corresponding author:} Statistics and
		Econometric Methods, Goethe University Frankfurt, Theodor-W.-Adorno-Platz 4,
		60323 Frankfurt, Germany, email: \texttt{hassler@wiwi.uni-frankfurt.de} .} \\
	Goethe University Frankfurt \and Mehdi Hosseinkouchack \\
	University of Mannheim }
\maketitle
\begin{abstract}
	\noindent The number $\frac{\pi ^{2}}{6}$ is involved in the variance of several distributions in statistics. At the same time it holds $\sum\nolimits_{k=1}^{\infty }k^{-2}= \frac{\pi ^{2}}{6}$, which solves the famous Basel problem. We first provide a historical perspective on  the Basel problem, and second show how to generate further proofs building on stochastic processes.
\end{abstract}



\section{Introduction}

The number $\frac{\pi ^{2}}{6}$  solves the so-called Basel problem in that $\sum\nolimits_{k=1}^{\infty }k^{-2}= \frac{\pi ^{2}}{6}$; the first proof was given by Leonhard Euler. A statistician encounters  this number  in several places.  First, $\frac{\pi ^{2}}{6}$ equals the variance of the famous standard Gumbel (or extreme value) distribution, see \citet{Gumbel1941}. Second, it shows up in the limiting variance of the appropriately normalized estimator from the so-called log-periodogram regression. This estimator of the fractional order of integration is sometimes called GPH estimator after the seminal paper by \citet{GPH83}, which was recently celebrated by a special issue of the \textit{Journal of Time Series Analysis}, see \citet{NielsenHualde2019}; confer also \citet*[Thm. 1]{HDB98}. Third, the inverse of this number, $\frac{6}{\pi ^{2}}$, amounts to the  variance in the limiting distribution of the maximum likelihood estimator of the order of integration of so-called fractionally integrated noise, see \citet[Coro. 8.1]{Hassler19}.

In this note, we show how to generate further proofs of $\sum\nolimits_{k=1}^{\infty }k^{-2}= \frac{\pi ^{2}}{6}$ using Karhunen-Lo\`{e}ve expansions of the Wiener process and of related processes. The next section briefly reviews the history of the Basel problem until the year 1735. Section 3 provides a short review and a classification of the multitude of earlier proofs. Section 4 uses a technique from the theory of stochastic processes to produce new proofs. The final section provides some concluding remarks.

\section{Basel problem}

According to \citet{enestroem1911}, it was Pietro Mengoli\footnote{Born in 1625 according to \citet{boyer1968}, while \citet{enestroem1911} dates his year of birth to 1626. In 1650 he published \emph{Novae Quadraturae Arithmetica, sev De Additione Fractionum}, which seems to have been largely forgotten until the rediscovery by Enestr\"{o}m.}  (1625/26-1686) of Bologna in Italy who posed the problem to determine the value of $\sum\nolimits_{k=1}^{\infty }k^{-2}$, for which we nowadays write $\zeta (2)$ with Riemann zeta function $\zeta(s)=\sum\nolimits_{k=1}^{\infty }k^{-s}$ for complex $s$ with $\mbox{Re} (s) >1$. It was clear that $\zeta (2)$ is finite for the following reason. By partial fractions the series of reciprocals of triangular numbers ``telescopes'',
\[
T_N := \sum_{k=1}^N \frac{2}{k (k+1)} = 2 \sum_{k=1}^N  \left( \frac{1}{k } - \frac{1}{k+1}\right) = 2 \, \left(1 - \frac{1}{N+1}\right) \to 2
\]
as $N \to \infty$. Because of $2 \, k^2 \geq k^2 + k$ it follows that $\zeta (2) \leq 2$. In fact, \citet[p. 144-145]{enestroem1911} gave three different proofs for $T_N \to 2$ originally published by Mengoli in 1650. It was Jakob Bernoulli (1654-1705) from the city of Basel in Switzerland who popularized such results (\emph{Tractatus de seriebus infinitis}, 1689), apparently without being aware of the previous work by Mengoli. We quote his treatise as \citet{bernoulli1713} published posthumously  together with his famous \emph{Ars Conjectandi}. His \emph{tractatus} is organized in 60 propositions.  \citet[Prop. XV]{bernoulli1713} established convergence of $\sum_{k=1}^N \frac{2}{k (k+1)}$, and he proved in Prop. XVI that the harmonic series $\sum\nolimits_{k=1}^{N}k^{-1}$ does not converges, but he failed in Prop. XVII to  determine the value of $\zeta (2)$; \citet[p. 254]{bernoulli1713} wrote the request that he would be much obliged if someone found what had escaped his efforts and communicated it to him. Hence,  the Basel problem was there and resisted the efforts of  mathematicians of this time. In 1705, Johann Bernoulli (1667-1748) took the position of his late brother Jakob at the University of Basel, see \citet{merian1860} for a history of the mathematical Bernoulli dynasty. Leonhard Euler (1707-1783) graduated at the University of Basel in 1724 at the age of 17 as one of Johann Bernoulli's students, see the biographical sketch in \citet[p. 1068]{ayoub1974}.  The Basel problem remained unsolved until \citet{euler1735} provided the following formula.\footnote{It was published in Volume 7 of the Memoirs of the Imperial Academy of Sciences in St. Petersburg (Commentarii Academiae Scientiarum Imperialis Petropolitanae). The title page of this volume, which by the way contains 10 papers by Euler, carries three year dates: 1734/1735 and 1740. The first ones are the years when the papers were presented to the Academy, and the latter date refers to  the year of printing; we follow the international custom and quote Euler's paper as \citet{euler1735}.}

\bigskip

\noindent \textbf{Theorem} $\zeta \left( 2\right) =\frac{\pi ^{2}}{6}$.

\bigskip

\noindent We briefly review a selection of earlier proofs next.

\section{Earlier proofs}

Euler's first and famous proof built on a factorization of the Taylor expansion of the sine function, see e.g. the exposition by \citet[Sect. 4]{ayoub1974}. According to \citet[p. 1077]{ayoub1974}, Daniel Bernoulli (1700-1782, son of Johann Bernoulli, see again \citet{merian1860}) did not doubt Euler's result, but had two concerns with respect to the validity of the arguments leading there. As early as 1738, Nicolaus Bernoulli (1687-1759, a nephew of Johann and Jacob Bernoulli) presented quite a different idea of proof, which  nowadays is rarely associated with him.\footnote{Just like Euler's original paper it was published in  the Memoirs of the Imperial Academy of Sciences in St. Petersburg, Volume 10. This time the title page  carries two year dates, 1738 and 1747, where 1747 is  the year of printing; we  refer to the paper as \citet{bernoulli1738}. We do not correct for a misprint in the title: Instead of  $\frac{1}{6}$ one should read of course  $\frac{1}{9}$.} Bernoulli (1738) defined the sum over reciprocals of odd squares as $Z_N$, and the alternating series $q_N$:
\[
Z_N = \sum_{k=0}^N \frac{1}{(2k+1)^2} \quad \mbox{and} \quad q_N = \sum_{k=0}^N \frac{(-1)^{k}}{(2k+1)} \, .
\]
Bernoulli's  idea was to square $q_N$, such that $Z_N = q_N^2 + y_N$, where $y_N$ is implicitly defined. Then he argued that $q_N \to \pi/4$ and $y_N \to \pi^2/16$  as $N \to \infty$, which amounts to
\begin{equation} \label{odd_squares}
\sum_{k=0}^\infty \frac{1}{(2k+1)^2} = \frac{\pi ^{2}}{8} \, .
\end{equation}
Obviously, this solves the Basel problem since
\begin{equation} \label{Euler2}
\zeta(2) - \frac{1}{4} \, \zeta (2) = \sum_{k=0}^\infty \frac{1}{(2k+1)^2}  \, .
\end{equation}
Similarly, with more rigorous arguments but without referring to \citet{bernoulli1738}, \citet[p. 12]{estermann1947} established $Z_N = 2 \, q_N^2 + u_N$, $u_N \to 0$, which proves (\ref{odd_squares}) again;  see also \citet[p. 322-324]{knopp1954theory}.

Euler himself added a second proof in 1743 that builds on (\ref{odd_squares}) with (\ref{Euler2}), although it does not rely on squaring $q_N$ but rather on an expansion of the arcsine function. \citet{stackel1907vergessene} called this paper  a ``forgotten treatise by Leonhard Euler'' and provided a historical perspective  including a reprint of the French paper.\footnote{The title is \emph{D\'{e}monstration de la somme the cette Suite ${1} + \frac{1}{4} + \frac{1}{9} + \frac{1}{16} + \frac{1}{25} + \frac{1}{36} + \mbox{etc.}$} published in Volume 2 of a journal  called somehow surprisingly \emph{Journal litt\'{e}raire d'Allemagne, de Suisse et du Nord}.} For Euler's second proof in a nutshell we refer to \citet{kimble1987euler} or footnote 63 in \citet[p. 376]{knopp1954theory}; see also  \citet{choe1987elementary}, where the last author does not seem to be  aware of reproducing Euler's second proof.

Ever since Euler cracked the Basel problem, many other proofs appeared. Some of them also culminate in establishing (\ref{odd_squares}), e.g. in chronological order: \cite{giesy1972still},  \cite{Hofbauer2002}, \cite{Harper2003}, \cite{ivan2008simple},  \cite{TimothyMarshall2010},   \cite{hirschhorn2011simple}, \cite{muzaffar2013new} and \cite{ritelli2013}. More proofs that do not establish (\ref{odd_squares}) in order to show $\zeta \left( 2\right) =\frac{\pi ^{2}}{6}$ were given by   \cite{knopp1918herleitung},\footnote{This proof can also be found in English in \citet[p. 266, 267]{knopp1954theory}.}  \cite{yaglom1953elementary}, \cite{matsuoka1961elementary}, \cite{stark1969another}, 
 and \cite{papadimitriou1973simple}. A note following \cite{papadimitriou1973simple} says, that this paper was translated from a  Greek manuscript, and that  the proof coincides with the one given in Norwegian by \cite{holme1970enkel};  we checked that it is actually identical with the one by  \cite{yaglom1953elementary}.\footnote{For simplicity, we quote \cite{yaglom1953elementary} with an English translation of the title although the paper is in Russian.}  Further, more or less elementary proofs have been published by \cite{apostol1983proof}, 
\cite*{beukers1993sums}, \cite{kortram1996simple}, \citet*{borwein1997pi},   \cite{passare2008compute},   \cite{daners2012short},\footnote{This proof is similar to the one by \cite{matsuoka1961elementary}, but 
more straightforward.}  
 \cite{xu2013connection} and \citet{lord2016most}. This list  does not include proofs that rely on Fourier analysis or proofs
that establish more generally closed form results for $\zeta (2n)$, $n \in \mathbb{N}$.  
 Therefore the list is far from being complete.

There are two recent proofs that stand out against all other ones  using a probabilistic approach: \cite{pace2011probabilistically} and \cite{holst2013probabilistic}. In the next section
we give additional proofs for $\zeta (2)=\frac{\pi ^{2}}{6}$ that are rooted in probability theory, too,  relying on  results from the theory of stochastic processes.

\section{Karhunen-Lo\`{e}ve expansions and new proofs} 

Consider a stochastic  zero mean process $X\left( t\right)$, $t \in [0,1]$, and assume that it is Gaussian and continuous in
quadratic mean with positive definite covariance kernel $k\left( s,t\right) = \mbox{E} (X(s) X(t)) $ that  is
symmetric and continuous. Then $X\left( t\right) $ is endowed with a Karhunen-Lo\`{e}ve\ (KL)
expansion,
\begin{equation} \label{KL}
X\left( t\right) =\sum_{j=1}^{\infty }\lambda _{j}^{-1/2}f_{j}\left(
t\right) Z_{j}\text{,}
\end{equation}%
where $Z_{j}$, $j\in 
\mathbb{N}$, is a sequence of independent standard normal random variables,  and  $\lambda _{j}$ and $f_{j}\left( t\right) $ are the eigenvalues and
eigenfunctions of $k\left( s,t\right) $ satisfying the following 
Fredholm integral equation, 
see \cite{loeve78} for a discussion: 
\begin{equation*}
f(t)=\lambda \int_0^1 k(s,t)f(s) \mbox{d} s.
\end{equation*}
It is worth noting that $f_{j}\left(
t\right) $, $j\in 
\mathbb{N}$, form an orthonormal base for $L^{2}$.  \cite{ito1968} established that the sum\ $\sum\nolimits_{j=1}^{N}\lambda
_{j}^{-1/2}f_{j}\left( t\right) Z_{j}$ converges a.s. uniformly to $X\left(
t\right) $, $N\rightarrow \infty$. From (\ref{KL})  we have 
\begin{equation} \label{KL_var}
k(t,t) =\sum_{j=1}^{\infty } \frac{f_{j}^2 (t)}{\lambda _{j}} \, ,
\end{equation}
see also Mercer's Theorem, e.\,g. \citet[Thm. 5.2]{Tanaka96}. For three related processes we will use (\ref{KL_var}) to show $\zeta \left( 2\right) =\frac{\pi ^{2}}{6}$. Proofs that are similar in spirit can be established using \citet[Coro. 2, p. 92]{Hochstadt73}.

Let $W\left( t\right) $ denote a standard Brownian motion or  Wiener process with kernel $k(s,t) = \min (s,t)$. Further, $W^\mu$ and $W^\tau$ are demeaned and detrended Wiener processes, respectively; they are defined as the orthogonal component of the projection on a constant or on a linear time trend:
\[
W^\mu (t) = W(t) - \int_0^1 W(s) \mbox{d} s \, ,
\]
\[
W^\tau (t) = W(t) + (6t-4) \int_0^1 W(s) \mbox{d} s + (6 - 12t) \int_0^1 s W(s) \mbox{d} s \, .
\]
For their kernels it is straightforward to obtain 
\[
k^\mu (s,t) = \min (s,t) - (s+t) + \frac{1}{2} (s^2+t^2)  + \frac{1}{3}  \, , 
\]
\[
 k^\tau (s,t) = \min (s,t)  - \frac{11}{10}(s+t) + {2} (s^2+t^2) - (s^3+t^3) -3 (st^2 +ts^2) +2 (st^3+ts^3) + \frac{6}{5} st +  \frac{2}{15} \, ,
\] 
where the latter result can be found e.g. in \citet*[Lemma  2.1]{ai2012karhunen}. 
Further, the eigenstructure  has been characterized as follows.

\bigskip

\noindent \textbf{Lemma} \textit{The eigenvalues and eigenfunctions of $k$, $k^\mu$ and $k^\tau$ are} ($j=1,2, \ldots$)
\[
\lambda_j = (j-1/2)^2 \pi^2 \, , \quad \lambda_j^\mu = j^2 \pi^2 \, , \quad \lambda_j^\tau = \left\{ \begin{array}{cc} (j+1)^2 \pi^2 \, , & j=2n-1 \\
4 z_{3/2, j/2}^2 \, , & j=2n \end{array} \right. ,
\]
\textit{where} $z_{3/2,n}$ \textit{are the positive roots of the Bessel function} $J_{3/2}$ \textit{of the first kind, and}
\[
f_{j}\left( t\right) =\sqrt{2}\sin ((j-1/2)\pi t) \, , \quad f_{j}^{\mu }\left( t\right) =\sqrt{2}\cos (j\pi t) \, ,
\]
\[
 f_{j}^{\tau }\left( t\right) =   \left\{ \begin{array}{cc} \sqrt{2}\cos ((j+1)\pi t) \, , & j=2n-1 \\ (-1)^{(j+2)/2}
 \sqrt{\Lambda_j} \sin(\sqrt{\lambda_j}(t-1/2))\, , & j=2n \end{array} \right. ,
\]
\textit{where} $\Lambda_j$ \textit{is given in \citet[eq. (22)]{HosseinkouchackHassler2016} and can be reduced to} $\Lambda_j = 2 \left(\sin \frac{\sqrt{\lambda_j}}{2}\right)^{-2}$.

\bigskip

\noindent \textbf{Proof} \citet{HosseinkouchackHassler2016} derive the eigenstructure of the kernels of demeaned and detrended Ornstein-Uhlenbeck processes; $W^\mu$ and $W^\tau$ are embedded as special cases, see \citet[Remark 1 and 2]{HosseinkouchackHassler2016}. In particular, $\lambda^\mu$ and $f^\mu$ are also given by \citet*[p. 2495]{beghin2005}, and  \citet[Thm. 1]{ai2012karhunen} derive $\lambda^\tau$, however, without giving  $f^\tau$. The case of the standard Wiener process is textbook knowledge. \quad $\blacksquare$

\bigskip

\noindent With the Lemma at hand, new proofs of $\zeta \left( 2\right) =\frac{\pi ^{2}}{6}$ are obvious.\\
Proof 1:  From (\ref{KL_var}) we obtain for $k(t,t)=t$ with $\lambda_j$ and $f_j$ that
\begin{equation*}
t=8\sum_{j=1}^{\infty }\frac{\sin ^{2}\left[ \left( j-1/2\right) \pi t\right]
}{\left( 2j-1\right) ^{2}\pi ^{2}} \, .
\end{equation*}%
Evaluation for $t=1$ amounts to (\ref{odd_squares}), which completes the proof by (\ref{Euler2}).\\
Proof 2: From (\ref{KL_var}) we obtain for $k^\mu(1,1)$ by the Lemma that 
 $\frac{1}{3}=\frac{2}{\pi^2} \sum_{j=1}^{\infty }\frac{1}{j^2}$, which is the required result.\\
 Proof 3: Note from the Lemma that $f_j^\tau (1/2) = 0$ for even $j$.  From (\ref{KL_var}) it hence follows for $k^\tau(1/2,1/2)$ that $\frac{1}{12}=\frac{2}{\pi^2} \sum_{n=1}^{\infty }\frac{1}{4 n^2} = \frac{\zeta (2)}{2 \pi^2}$, which proves the result.

\section{Concluding remarks} 

More proofs of the Basel problem can be produced following the route tackled here. All one needs is the expansion from (\ref{KL_var}) for some stochastic process. We employed expansions for the Wiener process, for  the demeaned Wiener process and for the detrended Wiener process.  As a further example, one may consider the  so-called Brownian bridge, for which 
the  required expansion can be found in \citet[pp. 213-214]{shorack2009empirical}.

\bibliographystyle{chicago}
\bibliography{Basel00}

\begin{thebibliography}{}

\bibitem[\protect\citeauthoryear{Ai, Li, and Liu}{Ai
  et~al.}{2012}]{ai2012karhunen}
Ai, X., W.~V. Li, and G.~Liu (2012).
\newblock Karhunen-{L}o\`{e}ve expansions for the detrended {B}rownian motion.
\newblock {\em Statistics \& Probability Letters\/}~{\em 82\/}(7), 1235--1241.

\bibitem[\protect\citeauthoryear{Apostol}{Apostol}{1983}]{apostol1983proof}
Apostol, T.~M. (1983).
\newblock A proof that {E}uler missed: evaluating $\zeta(2)$ the easy way.
\newblock {\em The Mathematical Intelligencer\/}~{\em 5\/}(3), 59--60.

\bibitem[\protect\citeauthoryear{Ayoub}{Ayoub}{1974}]{ayoub1974}
Ayoub, R. (1974).
\newblock Euler and the zeta function.
\newblock {\em The American Mathematical Monthly\/}~{\em 81\/}(10), 1067--1086.

\bibitem[\protect\citeauthoryear{Beghin, Nikitin, and Orsingher}{Beghin
  et~al.}{2005}]{beghin2005}
Beghin, L., Y.~Nikitin, and E.~Orsingher (2005).
\newblock Exact small ball constants for some {G}aussian processes under the
  {L}2-norm.
\newblock {\em Journal of Mathematical Sciences\/}~{\em 128\/}(1), 2493--2502.

\bibitem[\protect\citeauthoryear{Bernoulli}{Bernoulli}{1713}]{bernoulli1713}
Bernoulli, J. (1713).
\newblock {\em Ars conjectandi, opus posthumum. {A}ccedit {T}ractatus de
  seriebus infinitis, et epistola gallicé scripta de ludo pilae reticularis}.
\newblock Thurneysen.

\bibitem[\protect\citeauthoryear{Bernoulli}{Bernoulli}{1738}]{bernoulli1738}
Bernoulli, N. (1738).
\newblock Inquisitio in summam seriei $\frac{1}{1} + \frac{1}{4} + \frac{1}{6}
  + \frac{1}{16} + \frac{1}{25} + \frac{1}{36} + \mbox{etc}$.
\newblock {\em Commentarii Academiae Scientiarum Imperialis
  Petropolitanae\/}~{\em 10}, 19--21.

\bibitem[\protect\citeauthoryear{Beukers, Kolk, and Calabi}{Beukers
  et~al.}{1993}]{beukers1993sums}
Beukers, F., J.~A.~C. Kolk, and E.~Calabi (1993).
\newblock Sums of generalized harmonic series and volumes.
\newblock {\em Nieuw Archief voor Wiskunde (4)\/}~{\em 11\/}(3), 217--224.

\bibitem[\protect\citeauthoryear{Borwein, Borwein, and Dilcher}{Borwein
  et~al.}{1989}]{borwein1997pi}
Borwein, J.~M., P.~B. Borwein, and K.~Dilcher (1989).
\newblock Pi, {E}uler numbers, and asymptotic expansions.
\newblock {\em The American Mathematical Monthly\/}~{\em 96\/}(8), 681--687.

\bibitem[\protect\citeauthoryear{Boyer}{Boyer}{1968}]{boyer1968}
Boyer, C.~B. (1968).
\newblock {\em A History of Mathematics}.
\newblock Wiley.

\bibitem[\protect\citeauthoryear{Choe}{Choe}{1987}]{choe1987elementary}
Choe, B.~R. (1987).
\newblock An elementary proof of $\sum_{n= 1}^{\infty}
  \frac{1}{n^2}=\frac{\pi^2}{6}$.
\newblock {\em The American Mathematical Monthly\/}~{\em 94\/}(7), 662--663.

\bibitem[\protect\citeauthoryear{Daners}{Daners}{2012}]{daners2012short}
Daners, D. (2012).
\newblock A short elementary proof of $\sum \frac{1}{k^2}=\frac{\pi^2}{6}$.
\newblock {\em Mathematics Magazine\/}~{\em 85\/}(5), 361--364.

\bibitem[\protect\citeauthoryear{Enestr{\"o}m}{Enestr{\"o}m}{1911}]{enestroem1911}
Enestr{\"o}m, G. (1911).
\newblock {Zur Geschichte der unendlichen Reihen um die Mitte des siebzehnten
  Jahrhunderts}.
\newblock {\em Bibliotheca mathematica (3. Folge)\/}~{\em 12}, 135--148.

\bibitem[\protect\citeauthoryear{Estermann}{Estermann}{1947}]{estermann1947}
Estermann, T. (1947).
\newblock Elementary evaluation of $\zeta (2k)$.
\newblock {\em Journal of the London Mathematical Society\/}~{\em 22}, 10--13.

\bibitem[\protect\citeauthoryear{Euler}{Euler}{1735}]{euler1735}
Euler, L. (1735).
\newblock De summis serierum reciprocarum.
\newblock {\em Commentarii Academiae Scientiarum Imperialis
  Petropolitanae\/}~{\em 7}, 123--134.

\bibitem[\protect\citeauthoryear{Geweke and Porter-Hudak}{Geweke and
  Porter-Hudak}{1983}]{GPH83}
Geweke, J. and S.~Porter-Hudak (1983).
\newblock The estimation and application of long memory time series models.
\newblock {\em Journal of Time Series Analysis\/}~{\em 4}, 221--238.

\bibitem[\protect\citeauthoryear{Giesy}{Giesy}{1972}]{giesy1972still}
Giesy, D.~P. (1972).
\newblock Still another elementary proof that $\sum_{k=1}^\infty
  \frac{1}{k^2}=\frac{\pi^2}{6}$.
\newblock {\em Mathematics Magazine\/}~{\em 45\/}(3), 148--149.

\bibitem[\protect\citeauthoryear{Gumbel}{Gumbel}{1941}]{Gumbel1941}
Gumbel, E.~J. (1941).
\newblock The return period of flood flows.
\newblock {\em The Annals of Mathematical Statistics\/}~{\em 12}, 163--190.

\bibitem[\protect\citeauthoryear{Harper}{Harper}{2003}]{Harper2003}
Harper, J.~D. (2003).
\newblock Another simple proof of
  $1+\frac{1}{2^2}+\frac{1}{3^2}+...=\frac{\pi^2}{6}$.
\newblock {\em The American Mathematical Monthly\/}~{\em 110\/}(6), 540--541.

\bibitem[\protect\citeauthoryear{Hassler}{Hassler}{2019}]{Hassler19}
Hassler, U. (2019).
\newblock {\em Time Series Analysis with Long Memory in View}.
\newblock Wiley.

\bibitem[\protect\citeauthoryear{Hirschhorn}{Hirschhorn}{2011}]{hirschhorn2011simple}
Hirschhorn, M.~D. (2011).
\newblock A simple proof that $\zeta(2)=\frac{\pi^2}{6}$.
\newblock {\em The Mathematical Intelligencer\/}~{\em 33\/}(3), 81--82.

\bibitem[\protect\citeauthoryear{Hochstadt}{Hochstadt}{1973}]{Hochstadt73}
Hochstadt, H. (1973).
\newblock {\em Integral Equations}.
\newblock Wiley.

\bibitem[\protect\citeauthoryear{Hofbauer}{Hofbauer}{2002}]{Hofbauer2002}
Hofbauer, J. (2002).
\newblock A simple proof of $1+\frac{1}{2^2}+\frac{1}{3^2}+...=\frac{\pi^2}{6}$
  and related identities.
\newblock {\em The American Mathematical Monthly\/}~{\em 109}, 196--200.

\bibitem[\protect\citeauthoryear{Holme}{Holme}{1970}]{holme1970enkel}
Holme, F. (1970).
\newblock En enkel beregning av $\sum_{k= 1}^{\infty} \frac{1}{k^2}$.
\newblock {\em Nordisk Matematisk Tidskrift\/}~{\em 18}, 91--92.

\bibitem[\protect\citeauthoryear{Holst}{Holst}{2013}]{holst2013probabilistic}
Holst, L. (2013).
\newblock Probabilistic proofs of {E}uler identities.
\newblock {\em Journal of Applied Probability\/}~{\em 50\/}(4), 1206--1212.

\bibitem[\protect\citeauthoryear{Hosseinkouchack and Hassler}{Hosseinkouchack
  and Hassler}{2016}]{HosseinkouchackHassler2016}
Hosseinkouchack, M. and U.~Hassler (2016).
\newblock Powerful unit root tests free of nuisance parameters.
\newblock {\em Journal of Time Series Analysis\/}~{\em 37}, 533--554.

\bibitem[\protect\citeauthoryear{Hurvich, Deo, and Brodsky}{Hurvich
  et~al.}{1998}]{HDB98}
Hurvich, C.~M., R.~Deo, and J.~Brodsky (1998).
\newblock The mean squared error of {G}eweke and {P}orter-{H}udak's estimator
  of the memory parameter of a long-memory time series.
\newblock {\em Journal of Time Series Analysis\/}~{\em 19}, 19--46.

\bibitem[\protect\citeauthoryear{It{\^o} and Nisio}{It{\^o} and
  Nisio}{1968}]{ito1968}
It{\^o}, K. and M.~Nisio (1968).
\newblock On the convergence of sums of independent {B}anach space valued
  random variables.
\newblock {\em Osaka Journal of Mathematics\/}~{\em 5\/}(1), 35--48.

\bibitem[\protect\citeauthoryear{Ivan}{Ivan}{2008}]{ivan2008simple}
Ivan, M. (2008).
\newblock A simple solution to {B}asel problem.
\newblock {\em General Mathematics\/}~{\em 16}, 111--113.

\bibitem[\protect\citeauthoryear{Kimble}{Kimble}{1987}]{kimble1987euler}
Kimble, G. (1987).
\newblock {E}uler's other proof.
\newblock {\em Mathematics Magazine\/}~{\em 60\/}(5), 282--282.

\bibitem[\protect\citeauthoryear{Knopp}{Knopp}{1951}]{knopp1954theory}
Knopp, K. (1951).
\newblock {\em Theory and application of infinite series}.
\newblock {Blackie \& Son Limited}.

\bibitem[\protect\citeauthoryear{Knopp and Schur}{Knopp and
  Schur}{1918}]{knopp1918herleitung}
Knopp, K. and I.~Schur (1918).
\newblock {\"U}ber die {H}erleitung der {G}leichung $\sum_{n= 1}^{\infty}
  \frac{1}{n^2}=\frac{\pi^2}{6}$.
\newblock {\em Archiv der Mathematik und Physik, Ser. 3\/}~{\em 27}, 174--176.

\bibitem[\protect\citeauthoryear{Kortram}{Kortram}{1996}]{kortram1996simple}
Kortram, R. (1996).
\newblock Simple proofs for $\sum_{k= 1}^{\infty}
  \frac{1}{k^2}=\frac{\pi^2}{6}$ and $\sin x=x
  \prod_{k=1}^{\infty}(1-\frac{x^2}{k^2 \pi^2})$.
\newblock {\em Mathematics Magazine\/}~{\em 69\/}(2), 122--125.

\bibitem[\protect\citeauthoryear{Lo\`{e}ve}{Lo\`{e}ve}{1978}]{loeve78}
Lo\`{e}ve, M. (1978).
\newblock {\em Probability theory, vol. {II}}.
\newblock Springer.

\bibitem[\protect\citeauthoryear{Lord}{Lord}{2016}]{lord2016most}
Lord, N. (2016).
\newblock The most elementary proof that $\sum_{k= 1}^{\infty}
  \frac{1}{k^2}=\frac{\pi^2}{6}$.
\newblock {\em The Mathematical Gazette\/}~{\em 100\/}(549), 429--434.

\bibitem[\protect\citeauthoryear{Marshall}{Marshall}{2010}]{TimothyMarshall2010}
Marshall, T. (2010).
\newblock A short proof of $\zeta (2)= \frac{\pi^2}{6}$.
\newblock {\em The American Mathematical Monthly\/}~{\em 117\/}(4), 352--353.

\bibitem[\protect\citeauthoryear{Matsuoka}{Matsuoka}{1961}]{matsuoka1961elementary}
Matsuoka, Y. (1961).
\newblock An elementary proof of the formula $\sum_{k= 1}^{\infty}
  \frac{1}{k^2}=\frac{\pi^2}{6}$.
\newblock {\em The American Mathematical Monthly\/}, 485--487.

\bibitem[\protect\citeauthoryear{Merian}{Merian}{1860}]{merian1860}
Merian, P. (1860).
\newblock {\em Die Mathematiker Bernoulli: Jubelschrift zur vierten
  S\"{a}cularfeier der Universit\"{a}t Basel}.
\newblock Schweighauser'sche Universit\"{a}ts-Buchdruckerei.

\bibitem[\protect\citeauthoryear{Muzaffar}{Muzaffar}{2013}]{muzaffar2013new}
Muzaffar, H.~B. (2013).
\newblock A new proof of a classical formula.
\newblock {\em The American Mathematical Monthly\/}~{\em 120\/}(4), 355--358.

\bibitem[\protect\citeauthoryear{Nielsen and Hualde}{Nielsen and
  Hualde}{2019}]{NielsenHualde2019}
Nielsen, M.~O. and J.~Hualde (2019).
\newblock Special issue of the \textit{Journal of Time Series Analysis} in
  honour of the 35th anniversary of the publication of {G}eweke and
  {P}orter-{H}udak (1983): Guest editors' introduction.
\newblock {\em Journal of Time Series Analysis\/}~{\em 40}, 386--387.

\bibitem[\protect\citeauthoryear{Pace}{Pace}{2011}]{pace2011probabilistically}
Pace, L. (2011).
\newblock Probabilistically proving that $\zeta (2)= \frac{\pi^2}{6}$.
\newblock {\em The American Mathematical Monthly\/}~{\em 118\/}(7), 641--643.

\bibitem[\protect\citeauthoryear{Papadimitriou}{Papadimitriou}{1973}]{papadimitriou1973simple}
Papadimitriou, I. (1973).
\newblock A simple proof of the formula $\sum_{k= 1}^{\infty}
  \frac{1}{k^2}=\frac{\pi^2}{6}$.
\newblock {\em The American Mathematical Monthly\/}~{\em 80\/}(4), 424--425.

\bibitem[\protect\citeauthoryear{Passare}{Passare}{2008}]{passare2008compute}
Passare, M. (2008).
\newblock How to compute $\sum1/n^2$ by solving triangles.
\newblock {\em The American Mathematical Monthly\/}~{\em 115}, 745--752.

\bibitem[\protect\citeauthoryear{Ritelli}{Ritelli}{2013}]{ritelli2013}
Ritelli, D. (2013).
\newblock Another proof of $\zeta (2)= \frac{\pi^2}{6}$ using double integrals.
\newblock {\em The American Mathematical Monthly\/}~{\em 120\/}(7), 642--645.

\bibitem[\protect\citeauthoryear{Shorack and Wellner}{Shorack and
  Wellner}{1986}]{shorack2009empirical}
Shorack, G.~R. and J.~A. Wellner (1986).
\newblock {\em Empirical processes with applications to statistics}.
\newblock Wiley.

\bibitem[\protect\citeauthoryear{St{\"a}ckel}{St{\"a}ckel}{1907}]{stackel1907vergessene}
St{\"a}ckel, P. (1907).
\newblock {Eine vergessene Abhandlung Leonhard Eulers {\"u}ber die Summe der
  reziproken Quadrate der nat{\"u}rlichen Zahlen}.
\newblock {\em Bibliotheca mathematica (3. Folge)\/}~{\em 8}, 37--60.

\bibitem[\protect\citeauthoryear{Stark}{Stark}{1969}]{stark1969another}
Stark, E. (1969).
\newblock Another proof of the formula $\sum_{k= 1}^{\infty}
  \frac{1}{k^2}=\frac{\pi^2}{6}$.
\newblock {\em The American Mathematical Monthly\/}~{\em 76\/}(5), 552--553.

\bibitem[\protect\citeauthoryear{Tanaka}{Tanaka}{1996}]{Tanaka96}
Tanaka, K. (1996).
\newblock {\em Time Series Analysis: Nonstationary and Noninvertible
  Distribution Theory}.
\newblock Wiley.

\bibitem[\protect\citeauthoryear{Xu and Zhou}{Xu and
  Zhou}{2014}]{xu2013connection}
Xu, H. and J.~Zhou (2014).
\newblock The connection between the {B}asel problem and a special integral.
\newblock {\em Applied Mathematics\/}~{\em 5}, 2570--2584.

\bibitem[\protect\citeauthoryear{Yaglom and Yaglom}{Yaglom and
  Yaglom}{1953}]{yaglom1953elementary}
Yaglom, A.~M. and I.~M. Yaglom (1953).
\newblock An elementary derivation of the formulas of {W}allis, {L}eibniz and
  {E}uler for the number $\pi$.
\newblock {\em Uspekhi Matematicheskikh Nauk\/}~{\em 8\/}(5), 181--187.

\end{thebibliography}

\end{document}